\documentclass[a4paper,11pt]{amsart}
\usepackage{amssymb,amsfonts,amsxtra, mathrsfs,placeins,graphicx,verbatim}
\usepackage[all]{xy}
\xyoption{line}
\usepackage{fullpage}
\usepackage{euscript}
\newtheorem{theorem}{Theorem}[section]
\newtheorem{cor}[theorem]{Corollary}
\newtheorem{lem}[theorem]{Lemma}
\newtheorem{prop}[theorem]{Proposition}

\theoremstyle{definition}

\newtheorem{example}[theorem]{Example}
\newtheorem{defi}[theorem]{Definition}
\newtheorem{rem}[theorem]{Remark}

\numberwithin{equation}{section}

\DeclareMathOperator{\Hom}{Hom}
\DeclareMathOperator{\Aut}{Aut}
\DeclareMathOperator{\Ve}{Vert}

\DeclareMathOperator{\Edge}{Edge}
\DeclareMathOperator{\im}{Im}
\DeclareMathOperator{\Ker}{Ker}
\DeclareMathOperator{\Det}{Det}

\def\g{\mathfrak g}
\def\MC{\mathscr {MC}}
\newcommand{\V}{{\EuScript V}}

\def\G{\mathcal G}
\def\E{\mathcal E}

\def\P{\mathcal P}
\def\Q{\mathcal Q}
\def\B{\mathsf B}
\def\De{\Det^{d}}

\def\T{\mathbb T}

\def\In{\operatorname{In}}

\def\ra{\rightarrow}

\def\Z{\mathbb{Z}}
\def\S{\mathbb{S}}

\def\ground{\mathbf k}

\newcommand{\noproof}{\begin{flushright}\ensuremath{\square}\end{flushright}}

\newcommand{\Ass}{\mathcal{A}\!\mathit{ss}}
\newcommand{\Comm}{\mathcal{C}\!\mathit{om}}
\newcommand{\Lie}{\mathcal{L}\!\mathit{ie}}

\newcommand{\BV}{\mathsf{BV}}
\newcommand{\bigBV}{\widehat\BV}
\newcommand{\bv}{\mathsf{bv}}
\newcommand{\Flag}{\operatorname{Flag}}

\def\C{\mathcal C}
\def\O{\mathcal O}
\def\E{\mathcal E}
\def\F{\mathsf F}

\def\ev{\operatorname{ev}}

\def\id{\operatorname{id}}

\thanks{}

\begin{document}

\title[Minimal models...]{Feynman diagrams and minimal models for operadic algebras}
\author{J. Chuang  \and A.~Lazarev}
\thanks{The first author is supported by an EPSRC advanced research fellowship}
\address{Centre for Mathematical Science\\City University\\London EC1V 0HB\\UK}
\email{J.Chuang@city.ac.uk}
\address{University of Leicester\\ Department of Mathematics\\Leicester LE1 7RH, UK.}
\email{al179@leicester.ac.uk}
\begin{abstract}
We construct an explicit minimal model for an algebra over the
cobar-construction of a differential graded operad. The structure maps of this minimal model are expressed in terms of sums over decorated trees.  We introduce the appropriate notion of a homotopy equivalence of operadic algebras and  show that our minimal model is homotopy equivalent to the original algebra. All this generalizes and gives a conceptual explanation of well-known results for $A_\infty$-algebras. Further, we show that these results carry over to the case of algebras over modular operads; the sums over trees get replaced by sums over general Feynman graphs. As a by-product of our work we prove gauge-independence of Kontsevich's `dual construction' producing graph cohomology classes from  contractible differential graded Frobenius algebras.
\end{abstract}

\maketitle
\section{introduction}
The existence of a minimal model for an $A_\infty$-algebra was proved by Kadeishvili \cite{Kad}; the
precursor
of this result
 is the Sullivan theory of minimal models for rational homotopy types \cite{Su}, the latter being minimal models for $L_\infty$-algebras of a special type. More recently minimal models have
found applications in theoretical physics, for example in string field theory and quiver gauge theory, cf. e.g. \cite{Laz, To, AK, AF, Kajiura}.

In the later developments  it is important to have explicit formulas for the structure maps of minimal models. The formulas appearing in \cite{Mer}, \cite{markl} and \cite{KS} express the structure maps of minimal $A_\infty$-algebras in terms of sums over trees similar to those arising in the perturbative expansions of path integrals. In these quoted references the formulas are established by direct calculation; one of the aims of the present paper is to give conceptual and combinatorics-free proofs. In doing so we discovered that the result holds in considerably greater generality, namely for algebras over cobar-constructions of dg operads (or, more generally, \emph{cofibrant} operads).

Encountering a sum over trees, one naturally seeks an
interpretation of the corresponding sum with arbitrary graphs replacing
trees.
 It turns out that sums over graphs do appear in formulas for minimal models of algebras over \emph{modular} operads. In some special cases, namely, for algebras over a modular operad that is the Feynman transform of the naive closure of a cyclic operad, the terms corresponding to graphs with nontrivial fundamental group vanish and we are left with trees only. In particular we recover formulas for minimal models of symplectic $A_\infty$-algebras, cf. \cite{Kajiura}.

Our main tool is the $\BV$-resolution of an operad, introduced in the modular context in \cite{CL}.
In the present paper we have chosen to focus instead on ordinary operads, to make the work accessible to
a larger audience.
 However we stress that
 our methods carry over in a straightforward fashion to the modular case; the corresponding results for modular operads are stated later on in the paper.

To any operad $\O$ we associate another operad $\BV\O$ and a quasi-isomorphism of operads $\BV\O\ra\O$ which admits a right inverse $\O\ra\BV\O$. An algebra over $\BV\O$ is the same as an $\O$-algebra together with a \emph{Hodge decomposition}. A Hodge decomposition of a complex is a decomposition of it into its homology and a contractible part, subject to some natural axioms. A Hodge decomposition is an instance of a strong homotopy retraction data (cf. for example \cite{La}).

The operad $\BV\O$
contains the canonical cofibrant resolution of $\O$ given by the double cobar-construction which we denote by $\bv\O$; the latter is essentially the linear version of the Boardman-Vogt tree complex \cite{BV}.  While not cofibrant, the operad $\BV\O$ is rather close to the cofibrant operad $\bv\O$; this fact allows one to associate a minimal model to a Hodge decomposition. In fact, we are using this property of $\BV\O$  (implicitly or explicitly) in almost all our constructions.

We also introduce the notion of a homotopy equivalence of algebras over operads or modular operads and show that (nonminimal) operadic algebras are homotopy equivalent to their minimal models. In the case of $A_\infty$, $L_\infty$ or $C_\infty$ algebras this notion reduces to the familiar one of an $\infty$-quasi-isomorphism. We believe that it is of independent interest; as an application we show that for a modular operad $\O$ the $\O$-graph cohomology classes arising from contractible $\O$-algebras via  Kontsevich's dual construction \cite{K1} do not depend on the choice of contracting homotopy (gauge independence).

In this  paper we work mostly in the category of  $\Z/2$-graded
 vector spaces (also known as super-vector spaces) over a field $\ground$ of characteristic zero.
However all our results (with obvious modifications) continue to hold in the $\Z$-graded context.
The adjective `differential graded' will mean `differential $\Z/2$-graded' and will be abbreviated as `dg'.  All of our unmarked tensors are understood to be taken over $\ground$.
 For a $\Z/2$-graded vector space $V=V_0\oplus V_1$ the symbol $\Pi V$ will denote the \emph{parity reversion} of $V$; thus $(\Pi V)_0=V_1$ while $(\Pi V)_1=V_0$.

\subsection{Minimal models for $A_\infty$-algebras}\label{Ainf}

Before diving into our general constructions, we wish to motivate our operadic approach.
So we begin by
recalling the explicit formulas for minimal models of
$A_\infty$-algebras given by Merkulov \cite{Mer}, interpreted as sums indexed over
planar trees by Kontsevich and Soibelman \cite{KS}.

Let $A$ be an $A_\infty$-algebra, i.e., a dg vector space equipped with
odd structure maps $$m_n:(\Pi A) ^{\otimes n}\rightarrow \Pi A, \quad n\geq 2,$$
subject to the identities
\begin{eqnarray}\label{Ainfconstraint}
\sum
_{i+j+k=n}
m_{i+1+k} (\id^{\otimes i}\otimes m_j \otimes \id^{\otimes k})=0, \quad n\geq 1,
\end{eqnarray}
where $m_1$ is the differential of $\Pi A$.

Let $B$ be another $A_\infty$-algebra.
An $A_\infty$-morphism $f:A\to B$ is a collection of even maps
$f_n:(\Pi A) ^{\otimes n}\rightarrow \Pi B, \quad n\geq 1$
intertwining the structure maps on $A$ and $B$ in an appropriate way.
We say that $f$ is an $A_\infty$-isomorphism (resp. $A_\infty$-quasi-isomorphism) if $f_1$ is an isomorphism (resp. quasi-isomorphism)
of dg spaces.

The structure map $m_2$ on $A$ induces an associative product on
the homology $H(A)$ of $A$. The minimal model theorem states that it is possible to extend this product to an $A_\infty$-structure on
$H(A)$ in such a way that there exists a $A_\infty$-quasi-isomorphism
$f:H(A)\to A$.

We now describe Merkulov's approach, in which both the $A_\infty$-structure on $H(A)$ and the map $f$ are constructed explicitly.
Choose a decomposition $A=W\oplus K$ of $A$ as a direct sum of sub dg spaces, such that $K$ is acyclic, together with a contracting homotopy $h:K\to K$ such that $h^2=0$. Later we shall call this a \emph{Hodge decomposition} of $A$. We are most
interested in \emph{canonical} Hodge decompositions, where the differential vanishes on $W$, so that $W$ is identified with $H(A)$; such
decompositions always exist. The chosen decomposition is not required to be compatible with the $A_\infty$-structure in any way.

\begin{figure}[h]
\[
\xymatrix
@C=3.5ex@R=2.50ex@M=0.3EX
{
&&&&& \ar[ddd]^t \\
&&& \ar[ddrr]^t \\
&&&&&&&& \ar[dlll]^t\\
&&&&& m_3 \ar[dddl]^s\\
&&& \ar[ddr]^t \\
 \ar[ddrr]^t && \ar[dd]^t\\
 &&&& m_2 \ar[dd]^s & \ar[ddl]^t \\
 && m_2 \ar[drr]^s &&&&& \ar[dlll]^t \\
 &&&& m_4 \ar[dd]^t \\
 \\
 &&&& \\
}
\]
\caption{Definition of $\tilde{m}_T$.}
\label{planartree}
\end{figure}
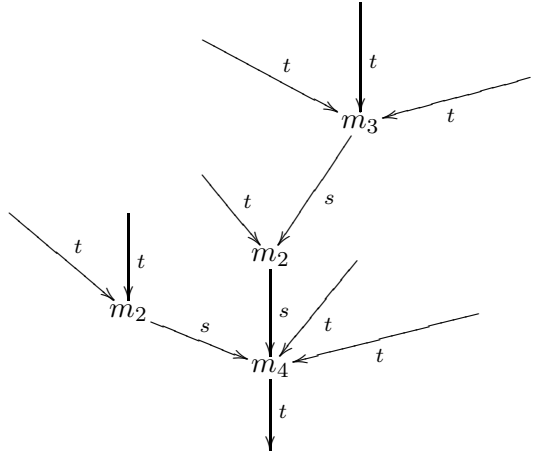

We define new operators
$\tilde{m}_n:(\Pi A) ^{\otimes n}\rightarrow \Pi A, \quad n\geq 2$
as follows. Let $t:\Pi A\to \Pi A$ be the projection of $\Pi A$
onto $\Pi W$ along $\Pi K$, and let $s:\Pi A\to \Pi A$ be equal to
$\Pi h$ on $\Pi K$ and $0$ on $\Pi W$. Let $T$ be a planar rooted tree with $n+1$ extremities; we assume that each vertex has valence at least $3$.
We label the extremities by $t$, all other edges by $s$ and each vertex of valence $v$ by $m_{v-1}$. Then
$\tilde{m}_T$ is constructed by working from the canopy of the tree down to the trunk, composing
labels in an obvious manner. For the tree pictured in Figure~\ref{planartree}
we have
\[
\tilde{m}_T=
tm_4
(\id \otimes sm_2 \otimes \id^{\otimes 2})
(sm_2\otimes \id\otimes sm_3 \otimes \id^{\otimes 2})
   t^{\otimes 8}.
\]

We define
$\tilde{m}_n = \sum_T \tilde{m}_T: (\Pi A)^{\otimes n}\to \Pi A$, where
the sum is taken over all planar rooted trees with $n+1$ extremities.
It can be shown by direct calculation \cite{Mer, markl, KS} that
the $\tilde{m}_n$ satisfy the $A_\infty$-constraint (\ref{Ainfconstraint}).
The new $A_\infty$-structure on $A$ thus defined restricts to one on $W$, and the inclusion of $W$ into $A$ extends to an $A_\infty$-quasi-isomorphism defined similarly as a sum over trees. In the case of canonical Hodge decomposition, $W\cong H(A)$, and the minimal model theorem is proved. (On the other hand for the trivial Hodge decomposition $A=A\oplus 0$, we obtain $\tilde{m}_n=m_n$.)

What is the significance of the individual operators $\tilde{m}_T$, before they are summed to obtain $\tilde{m}_n$ ? And what does it mean
to allow any edge of a tree to be labelled either by $s$ or by $t$ ?
In answering these questions, we will arrive
at the structure maps $\tilde{m}_n$ in a conceptual way, so that the $A_\infty$-constraint
(\ref{Ainfconstraint}) is automatically satisfied. A key idea is to regard the operators $s$ and $t$ coming from a Hodge decomposition of $A$ as part of an enhanced algebraic structure on $A$.

The language of operads is well suited to encode and develop the additional structure. We view an
$A_\infty$-algebra as an algebra over a certain operad $\O=\B\Ass$, and interpret a Hodge decomposition as an extension of the action
of $\O$ to a larger operad $\BV\O$.  The operations $\tilde{m}_T$ represent
the actions of particular elements of $\BV\O$; in this context sums over trees arise naturally. Moreover the passage from an $A_\infty$-algebra $A$ to its minimal model $H(A)$ may be regarded as a homotopy from the trivial Hodge decomposition to a canonical Hodge decomposition.

\subsection{Notation and conventions}

The general modern reference for differential graded operads is the work of Ginzburg and Kapranov \cite{GK}; the corresponding reference for \emph{modular} operads is \cite{GeK}. We adopt most of the notation and terminology from these two papers.

An $\S$-module is a collection of dg vector spaces $\{\V(n)\left|\right. n\geq 1\}$ with an action of the symmetric group $\S_n$ on $\V(n)$.
   Furthermore, if $\V$ is an $\S$-module and $I$ is a finite set then we set
\[\V(I):=\Bigl[\bigoplus\V(n)\Bigr]_{\S_n}\]
where the direct sum is extended over all bijections $\{1,2,\ldots,n\}\ra I$.

Recall that a \emph{dg operad} is an $\S$-module $\O$ together with composition maps
$\circ_i:\O(m)\otimes\O(n)\rightarrow \O(m+n-1)$ satisfying some natural
compatibility and invariance conditions. For brevity's sake, when
we say `operad' we will usually mean `dg operad'.

 An operad $\O$ is \emph{unital} if there is an element $\mathbf{1}\in\O(1)$ such that $\mathbf{1}\circ_1 x = x = x\circ_i \mathbf{1}$ for
all $x\in\O(m)$, $1\leq i \leq m$.
A unital operad $\O$ is called \emph{admissible} if $\O(1)=\ground$ is
the span of $\mathbf{1}$, and $\O(n)$ is a finite-dimensional dg vector space for any $n$. This notion is slightly more restrictive than the one used in \cite{GK} but sufficient for our purposes.

For an admissible operad $\O$ we will denote by $\bar\O$ the non-unital operad for which \[\bar\O(n)=\begin{cases} 0 & \text{if } n=1,\\ \O(n) &\text{if }n>1.\end{cases}\]

For a dg vector space $V$ we will write $\E(V)$ for the endomorphism operad $\{\E(V)(n)\}=\{\Hom(V^{\otimes n},V)\}$ of $V$.
An \emph{algebra structure} over a dg operad $\O$ on $V$ is a map of dg operads $\O\ra\E(V)$. If $\O$ is unital we usually assume the map is unital, i.e., the unit of $\O$ acts as the identity on $V$.

In this paper the language of \emph{trees} is used throughout. A tree is a non-empty oriented connected graph $T$ with no loops such that any vertex of $T$ has exactly one outgoing edge and at least one incoming edge. The edges which abut only one vertex are called \emph{extremities}; the unique outgoing extremity is called the \emph{root} of a tree and the remaining extremities are called \emph{leaves}. The set of \emph{internal edges}, i.e., edges which are
not extremities, is denoted $\Edge(T)$. The collection of vertices of $T$ will be denoted by $\Ve(T)$. For a vertex $v\in T$ the set of incoming edges will be called $\In(v)$. A tree is \emph{reduced} if it has no bivalent vertices.
Let $I$ be a finite set. An \emph{$I$-labelled tree} is a tree $T$ together with a bijection between $I$ and the set of leaves of $T$. A $\{1,\ldots,n\}$-labelled tree will simply be called an \emph{$n$-tree}.

If $\V=\{V(n)\}$ is an $\S$-module and $T$ is a tree we will denote by $\V(T)$ the dg vector space $\bigotimes_{v\in\Ve(T)}\V(\In(v))$; we will also call $\V(T)$ the \emph{space of $\V$-decorations on $T$} since it is spanned by tensors
$\otimes_{v} x_v$ corresponding to a choice of a `decoration' $x_v\in\V(\In(v))$
on each vertex $v$ of $T$.

The free operad on an $\S$-module $\V$ will be denoted by $\T\V$; recall that $\T\V(n)=\oplus_T \V(T)$, where the sum is over all (isomorphism classes)
of $n$-trees.

For any $n$-tree $T$ an operad $\O$ determines a homomorphism $\mu_T:\O(T)\ra\O(n)$ which corresponds to taking operadic compositions in $\O(T)$ along the internal edges of $T$.

For an ungraded vector space $V$ of dimension $n$
and $d\in \Z/2$
we write $\De(V)$ for the $\Z/2$-graded vector space $(\Pi^n\Lambda^n V)^{\otimes d}$.
 For a finite set $S$ we write $\De(S)$ for
$\De(\ground^S)$. Note that $\De(S)^*$ is canonically
isomorphic to $\De(S)$.

The paper is organized as follows. In Section \ref{recall} we recall the notions of the cobar-construction of an operad, as well as the $\bv$-- and $\BV$--resolutions of an operad. We also discuss the analogous constructions in the context of modular operads. Section \ref{main} contains our main construction -- a minimal model of an algebra over an admissible operad. In Section \ref{sechomotopy} we show that our minimal model is indeed homotopy equivalent (in the appropriate sense) to the original operadic algebra. This general notion of homotopy equivalence agrees with the familiar notion of infinity-equivalence in the case of $A_\infty, C_\infty$ or $L_\infty$-algebras. We also discuss the Maurer-Cartan moduli spaces associated to a differential graded Lie algebras and a suitable notion of homotopy between Maurer-Cartan elements. This material is presumably well known to experts but we are unaware of any published reference. In section \ref{modular} we extend our results to the setting of algebras over modular operads.

\section{Resolutions of operads}\label{recall}
We start by recalling the notion of the cobar-construction of an operad, following \cite{GK}.
\begin{defi}
Let $\P$ be an admissible dg operad. The cobar-construction $\B\P$ is the dg operad whose underlying operad of graded vector spaces is the free operad on the $\S$-module 
$\Pi\bar\P^*$.
 The differential in $\B\P$ is the sum of the internal differential in $\bar\P^*$ and the cobar-differential. The latter is induced on
 $\Pi\bar\P^*$
  by  the structure map $\T\P\rightarrow \P$ of the operad $\P$; it then uniquely extends to the whole of $\B\P$ by the Leibniz rule.
\end{defi}

Since $\bar\P(1)^*=0$, we have $\B\P(n)=\bigoplus_T (\Pi\bar\P^*)(T)$ where
the sum is over all \emph{reduced} $n$-trees. Moreover $\P\mapsto\B\P$ defines
a self-adjoint endofunctor on the category of admissible operads.
It is known, cf. \cite[Theorem (3.2.16)]{GK}, that the canonical counit map
$\B\B\P\rightarrow \P$ is a quasi-isomorphism, i.e.  $\B\B\P$ is a resolution of $\P$.

The structure of $\B\B\P$ appears rather complicated. We now describe another resolution of $\P$, the so-called $\BV$-resolution, which is very close to $\B\B\P$ but admits an extremely simple description in terms of generators and relations.

\begin{defi}Let $\P$ be an admissible dg operad. Its BV-resolution $\BV\P$ is the unital dg operad freely generated over $\P$ by an odd operation $s$ and an even operation $t$, both in $ \P(1)$, subject to the relations:\begin{itemize}\item
$s^2$=0;\item $t^2=t$;\item $st=ts=0$.
\end{itemize}
The differential $d$ on $\BV\P$ extends the differential on $\P$, and $d(s)=\mathbf{1}-t$ and  $d(t)=0$.
\end{defi}
\begin{rem}
We can write $\P[s,t]/(s^2,t^2-t,st)$ for the BV-resolution of $\P$ (without the differential).\end{rem}

Algebras over $\BV\P$ admit an especially simple description, as the following Proposition demonstrates; its proof is a simple unraveling of the definitions.

\begin{defi}\label{defHodge}
Let $V$ be a dg vector space. A \emph{Hodge decomposition} of $V$ is a choice of an odd operator $s:V\ra V$ such that
\begin{itemize}
\item
$s^2=0$,
\end{itemize}
and an even operator $t:V\ra V$ such that
\begin{itemize}
\item $t^2=t$;
\item $dt=td$.\end{itemize}
In addition the following identities for the operators $s$ and $t$ hold:
\begin{itemize}\item
$st=ts=0$;\item
$(ds+sd)(a)=a-t(a)$ for any $a\in V.$
\end{itemize}
If $dt=0$ we say that the Hodge decomposition is \emph{canonical}, and
if $t=\id_V$ (and thus $s=0$) that it is \emph{trivial}.
\end{defi}

\begin{prop}\label{dualalg1}
Let $\P$ be an admissible dg operad and $V$ be a dg vector space. Then the structure of an $\BV\P$-algebra on $V$ is equivalent to the following data:
\begin{enumerate}\item
The structure of a $\P$-algebra on $V$.
\item A Hodge decomposition of $V$.
\end{enumerate}
\end{prop}
\noproof

\begin{rem}\label{Hodge}
Given a Hodge decomposition of a dg vector space $V$, we have
$V=\im(t)\oplus\im(\operatorname{id}_V-t)$, and $s$ restricts to a contracting homotopy on the
second summand. Conversely any decomposition $V=W\oplus K$ where $K$ is equipped with a square-zero contracting homotopy determines a Hodge decomposition. The decomposition is canonical if and only if
$W\cong H(V)$ and trivial if and only if $K=0$.

Moreover canonical Hodge decompositions always exist.
Let $d$ denote the differential in $V$ and set  $V_0=\Ker d$ and $U=\im d$. Choose a complement $W$ to $U$ inside $V_0$, and a complement
$U'$ to $V_0$ inside $V$. We have therefore
\[V=W\oplus(U\oplus U^\prime).\]
Define the operator $t:V\ra V$ to be the projection onto $W$, and
define $s:V\ra V$ to be zero on $W\oplus U^\prime$ and inverse to $d: U^\prime\rightarrow U$ when restricted to $U$.
Then it is easy to check that $s$ and $t$ satisfy the identities of
Definition~\ref{defHodge}. Since $dt=0$ by construction, we obtain a canonical
Hodge decomposition.
\end{rem}

We will now give another, more concrete, description of the BV-resolution of an admissible dg operad as a kind of decorated tree complex. We start by defining the notion of a \emph{two-colored tree} or \emph{BV-tree}.

\begin{defi} A \emph{BV-tree} is a tree $T$ having the following additional structure:
the set $\Edge(T)$ of internal edges is partitioned into two subsets consisting of \emph{black edges} and \emph
{white edges}. We additionally require that the bivalent vertices of $T$ are
adjacent to extremities.
\end{defi}

We will denote the set of black edges of a BV-tree $T$ by $\Edge_b(T)$. For typographic reasons we will draw the black edges using straight lines and the white edges using wiggly lines. Note that the \emph{extremities} of a BV-tree have no color; to emphasize this we will draw the extremities using dotted lines.

There is an obvious notion of isomorphism between two BV-trees. Define two types of operations on BV-trees:\begin{enumerate}\item
contractions of black edges. For a black edge $e\in\Edge_b(T)$ this operation will be written as $T\mapsto T_e$.\item
replacing a black edge by a white edge. For a black edge $e\in\Edge_b(T)$ this operation will be written as $T\mapsto T^e$.\end{enumerate}

We can now introduce the two-colored tree version of the BV-resolution of an admissible operad $\P$, which will be temporarily  denoted by $\BV^\prime(\P)$.

 Given a BV-tree $T$ we put
 $$\P[T]:=\Det(\Edge_b(T))\otimes\P(T),$$
 a twisted version of the space of $\P$-decorations on $T$. To make sense of
 $\P(T)$ we are forgetting the coloring, treating $T$ as a usual tree.
  Let  $e\in \Edge_b(T)$ be a black edge. Then the contraction $T\mapsto T_e$ determines a parity-reversing linear map
$d_e:\P[T]\ra\P[T_e]$ given by the operadic composition in $\P$; similarly the operation $T\mapsto T^e$ determines tautologically an odd map $d^e:\P[T]\ra\P[T^e]$.

\begin{figure}[h]
\[
\xymatrix
@C=3.5ex@R=2.50ex@M=0.3EX
{
\ar@{.>}[dd] & \ar@{.>}[ddr] && \ar@{.>}[ddl] &&&& \\
&&&&&&&&&& \\
\bullet \ar[dr] && u \ar@{~>}[dl] & \ar@{.>}[ddr] && \ar@{.>}[ddl] & \\
& v \ar@{~>}[dd] &&&&&&&&
    &&\ar@{.>}[dd] & \ar@{.>}[ddr] && \ar@{.>}[ddl] \\
&&&& w \ar@{.>}[dd] &&&&&
    &&&&&&& \\
& \bullet \ar@{.>}[dd] &&&&&&&&
    &&\bullet \ar[dr] && u \ar@{~>}[dl] & \ar@{.>}[ddr] && \ar@{.>}[ddl]\\
&&&&&&&&&
    &&& v \ar@{~>}[dddrr] &&&&& \ar@{.>}[dll]\\
&&&& \circ &&&& =
    &&&&&&& x\circ w \ar[ddl] \\
& \circ &&& \ar@{.>}[dd] &&&&
    &&&&&& \\
& \ar@{.>}[dd] &&&&& \ar@{.>}[dll] &&&
    &&&&& y \ar@{.>}[dd] \\
&&&& x\ar[ddl] &&&&&
    &&&&&&& \\
& \bullet \ar@{~>}[drr] &&&&&&&&
    &&&&&&& \\
&&& y \ar@{.>}[dd] &&&&&& \\
&&&&&&&&&& \\
&&&&&&&&&& \\
}
\]

\caption{Composition of decorated BV-trees}
\label{BVcomposition}
\end{figure}
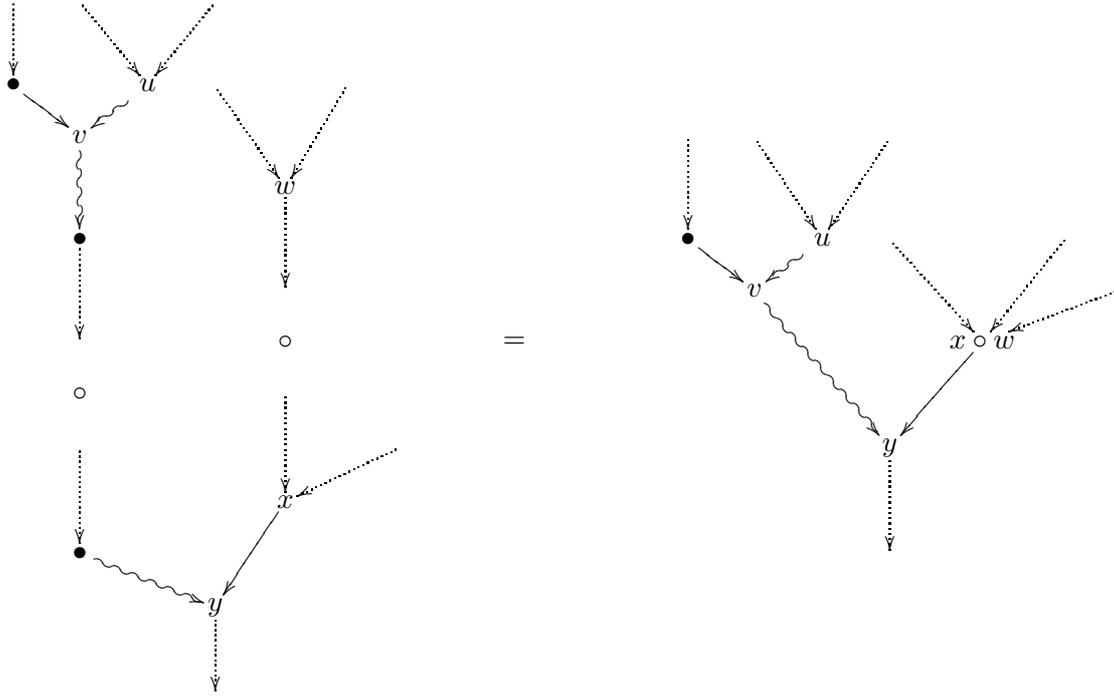

We define  \[\BV\P^\prime(n):=\bigoplus_{T} \P[T].\]
Here the direct sum is extended over isomorphism classes of all $n$-BV-trees $T$,
with differential $d$ determined by the formula
\begin{equation}\label{diff1}d|_{\P[T]}=d_\P+\sum_{e\in\Edge_b(T)}[d_e+d^e],\end{equation}
where $d_\P$ is the internal differential induced by the differential on $\P$. The following result is an analogue of Proposition 8.6 of \cite{CL}; its proof (which we omit) is similar to, but simpler than, the proof in the cited reference.
\begin{prop}\label{twoBVs}
There is an isomorphism of complexes
\[\BV\P(n)\cong\BV\P^\prime(n).\]
Under this isomorphism the operadic composition in $\BV\P$ is given by glueing decorated BV-trees according to the following rule: graft extremities to make a new internal edge and then contract the newly formed edge, using the operadic composition in $\P$. If the resulting decorated colored tree contains a bivalent vertex connecting two internal edges, it is considered to be zero, unless both edges are white in which case the vertex is removed and the edges are merged into a single white edge; this process ensures that a decorated BV-tree is obtained.
\end{prop}
\noproof
We  present in Figure~\ref{BVcomposition}
an example of a composition of three decorated BV-trees. The bivalent vertices are implicitly decorated by $\mathbf{1}\in\ground=\P(1)$.

We now explain how to compute
\emph{BV-tree amplitudes} of $\BV\P$-algebras. Let $V$ be a dg vector space with the structure of a $\BV\P$-algebra; this amounts to having a $P$-algebra structure on $V$ together with two operators $s$ and $t$ on $V$ satisfying the conditions of Definition~\ref{defHodge}. Let $x\in\P[T]$ be a $\P$-decoration on an $n$-BV-tree $T$; Then via the action map $\BV\P\to\E(V)$, $x$ determines an operator $Z^\P_V(x)\in\E(V)(n)=\Hom(V^{\otimes n}, V)$ with the following description.

The action map $\P\ra\E(V)$ allows one to view $x$ as an $\E(V)$-decoration on $T$. Thus, each vertex $v$ of $T$ has an element in $\E(V)(\In v)$ attached to it. Take the tensor product of these elements over all vertices of $T$. Then $Z^\P_V(x)$ is computed by contracting the resulting tensor along all edges of $T$ and interpreting black edges and white edges as compositions with the operators $s$ and $t$ respectively.


We conclude by explaining how $\BV\P$ is closely related to the standard free resolution of $\P$ provided by the twice-iterated cobar-construction.
To any BV-tree $T$  we can associate the BV-tree $T_t$ obtained from $T$ by glueing a white edge onto each leg. Then the image of the map $T\mapsto T_t$ consists of the BV-trees all of whose extremities are connected to white edges by bivalent vertices.

We define $\bv\P$ to be the truncation of $\BV\P$ by the idempotent $t$, i.e. the suboperad of $\BV\P$ defined by
\[\bv\P(n):=\bigoplus_T \P[T_t],\] where the sum is over all $n$-BV-trees.
The element $t$ serves as the operad unit for $\bv\P$.
The following result is the tree version of Proposition 8.8 of \cite{CL}; the proof carries over almost verbatim.
\begin{prop} Let $\P$ be an admissible dg operad.
\begin{enumerate}
\item
We have an isomorphism
$\bv\P\cong\B\B\P$ of dg operads.
\item The inclusion $\bv\P\hookrightarrow\BV\P$ is a quasi-isomorphism of dg operads.
\item
The embedding $\P\hookrightarrow\BV\P$ and the augmentation
$\BV\P\ra \P$ are quasi-isomorphisms of dg operads.
\end{enumerate}
\end{prop}
\noproof

\section{Main construction}\label{main}
Let $\O$ be an operad and $V$ be an $\O$-algebra; by abuse of notation we shall also denote by $V$ the underlying dg vector space. The following diagram of operads summarizes the results of the previous section:
\begin{equation}\label{operadmaps}\xymatrix {\bv\O\ar^i[r]\ar_p[d]&\BV\O\ar@{-->}^h[d]\ar^q@/^/[dl]\\ \O\ar_f[r]\ar^j@/^/[ur]&\E(V)&\E(H(V))\ar@{_{(}->}[l]}\end{equation}
The map $q:BV\O=\O[s,t]/(s^2,st,t^2-t)\rightarrow \O$ is defined by setting $s\mapsto 0, t\mapsto 1$; the map $j$ is the canonical splitting of $q$ given by considering elements of $\O$ as $\O$-decorated BV-trees with no edges. The map $i$ is the inclusion of $\bv\O$ as a suboperad consisting of $\O$-decorated trees whose extremities are connected to white edges by bivalent vertices, and $p=q\circ i$. The map $f$ is the given $\O$-algebra structure on $V$. The map $h$ corresponds to a chosen Hodge decomposition of $V$. Recall that the maps $i, j, q$ and $p$ are quasi-isomorphism of dg operads and that the map $p:\bv\O\rightarrow\O$ is the canonical resolution of $\O$ by the cofibrant operad $\bv\O$. Note also that $h\circ j=f$, and hence the whole diagram (\ref{operadmaps}) is homotopy commutative.

Finally, in the case when our Hodge decomposition is canonical the projector $t:V\rightarrow H(V)$ determines an inclusion of dg operads $\E(H(V))\hookrightarrow\E(V)$.

The composition $h\circ i$ allows one to regard $V$ as a $\bv\O$-algebra.
The action of $\bv\O$ restricts to the image $\im(t)$ of the operator $t$ of the Hodge decomposition, and if
the chosen Hodge decomposition of $V$ is canonical then $\im(t)=H(V)$. In this case it is natural to call the space $H(V)$ together with this $\bv\O$-algebra a \emph{minimal model} of the $\O$-algebra $V$. We will, however, reserve the term `minimal model' for another (closely related) notion.
\begin{defi}Suppose that the map of operads $p:\bv\O\rightarrow \O$ admits a right inverse, i.e., a map $k:\O\to\bv\O$ so that $p\circ k=\operatorname{id}_\O$, and that the Hodge decomposition of the dg vector space $V$ is canonical. The structure of an $\O$-algebra on the space $V$ determined by the operad map $h\circ i\circ k:\O\rightarrow \E(V)$ is called a minimal model of the $\O$-algebra $V$. Since the image $h\circ i\circ k:\O\rightarrow \E(V)$ is contained in $\E(H(V))\hookrightarrow \E(V)$  we will
use the term minimal model
also to refer to the corresponding $\O$-algebra structure on $H(V)=\im(t)\subset V$.
\end{defi}

We are most interested in the case when the operad $\O$ is the cobar-construction $\O=\B\P$ of an operad $\P$.  For example for $\P=\Ass, \Comm, \Lie$ the $\O$-algebras are (up to parity reversion) $A_\infty, L_\infty, C_\infty$-algebras respectively. Since $\B\P$ is a free operad on $\Pi
\overline{P}^*$ (disregarding the differential) the operad map $\B\P\rightarrow \E(H(V))$ specifying a minimal model of $V$ is determined by a collection of maps \[\Pi\P(n)^*\rightarrow \E(H(V))(n)=\Hom(H(V)^{\otimes n},H(V)),\quad n\geq 2.\]  These maps will be called the structure maps of the corresponding minimal model. For example if $\P=\Comm$ (so that $\P(n)=\ground$ is the trivial representation of
$\S_n$) we simply have a collection of $\Sigma_n$-equivariant odd maps $H(V)^{\otimes n}\ra H(V)$ determining the structure of a (minimal) $L_\infty$-algebra on $H(\Pi V)$.

To any reduced tree $T$ we associate a BV-tree $T_{\BV}$ by the following
recipe: color the edges of $T$ black and then glue a new white edge onto each leg. The map $T\mapsto T_{BV}$ is a bijection of the set of all reduced trees onto the set of BV-trees such that
\begin{enumerate}\item
each leg of $G$ abuts a bivalent vertex;\item
all edges of $G$ adjacent to extremities are white;\item
all other edges of $G$ are black.
\end{enumerate}

Recall that operadic composition determines a map $\mu_T:\P(T)\ra\P(n)$ for
each reduced $n$-tree $T$; denote by $\mu_T^*$ the $\ground$-linear  dual map.
We have an isomorphism $\Pi(\P(T)^*)\cong (\Pi \P^*)[T_{\BV}]$ and hence
a natural inclusion $\iota_T:\Pi(\P(T)^*) \hookrightarrow \B\P[T_{\BV}]$.

\begin{theorem}\label{main1}
Any choice of a canonical Hodge decomposition on a $\B\P$-algebra $V$  gives rise to a $\B\P$-algebra structure on $H(V)$ that is a minimal model of $V$.  The structure maps
\[
m_n:\Pi\P(n)^*\rightarrow\Hom(H(V)^{\otimes n},H(V))
\]
of this minimal model
 are given as follows:
\[
m_n=\sum_T  Z^{\B\P}_V \circ \iota_T \circ \mu^*_T,
\] where the summation is extended over all reduced $n$-trees $T$.
\end{theorem}
\begin{proof}
A canonical splitting $k:\B\P\rightarrow \B\B\B P$ is obtained by applying the contravariant functor $\B$ to the canonical projection (counit)
$p:\bv\P = \B\B\P\rightarrow \P$.

Interpreting $\B\B\B\P$ as a subspace of the space of $\B\P$-decorated BV-trees we see that the image of the restriction of $k$ to $\Pi\P^*$ involves only BV-trees satisfying conditions (1), (2)  and (3) above;
 moreover with this identification we have $k(x)=\mu^*(x)$ for any element $x\in\Pi\P^*\subset\B\P$. From this the formula for the structure maps is immediate.
\end{proof}

In other words, we obtain a minimal model in which the action of
$x\in\Pi\P^*$ is described as follows. Write $\mu^*(x)$ as a sum
of $\Pi\P^*$-decorated reduced trees $\mu^*_T(x)$. Label each leg of $T$ by
$t$ and each internal edge by $s$ and calculate
an amplitude as in Section~\ref{Ainf}, interpreting the decorations on
vertices as multilinear maps on $V$ via the originally specified
action of $\B\P$ on $V$. Summing these amplitudes over all $T$ with $n$ leaves we obtain the $n$-th component of the action of
$x$ on $H(V)$.

\begin{rem}
The statement of the theorem simplifies considerably for $\P=\Comm$ or $\P=\Ass$. Indeed, in these cases a $\P^*$-decorated tree is a tree or a planar tree respectively. For example, to obtain the structure map $\tilde{m}_n:(\Pi V)^{\otimes n}\ra \Pi V)$ of a minimal $A_\infty$-model one associates the (original) structure maps $m_i$'s to the corresponding vertices  of planar trees; the operator $t$ to the extremities, the operator $s$ to the internal edges and takes the sum of tree amplitudes over all planar trees with $n$ leaves. That recovers the formula derived in \cite{Mer}, \cite{KS} by a different method (see Section~\ref{Ainf}).
\end{rem}
\begin{rem}
Let $f:\B\P\rightarrow \E(V)$ be a map giving a dg vector space $V$ the structure of a $\B\P$-algebra; a minimal model is given by the map $f'=h\circ i\circ k:\B\P\rightarrow \E(V)$, where $k:\O\to\bv\O$ is a right inverse to
$p:\bv\O\rightarrow \O$.
 It follows from the homotopy commutativity of (\ref{operadmaps}) that the maps $f$ and $f^\prime$ are chain homotopic, and in particular they induce the same maps on homology. There is, however, a much stronger constraint on the maps $f$ and $f^\prime$ -- they are homotopic as \emph{operad maps}. The notion of homotopy for operads and the corresponding equivalence relation will be considered in the next section.
\end{rem}

\begin{rem}
One advantage of our conceptual approach to explicit minimal models is
its manifest functoriality. Let $f:\P\to\Q$ be a morphism of operads.
Any $\B\P$-algebra $V$ may be viewed as a $\B\Q$-algebra by restriction along
$\B f:\B\Q\to\B\P$. Fixing a Hodge decomposition of $V$ we obtain both
$\B\P$- and $\B\Q$-algebra structures on $H(V)$, and by construction
the two are again related by restriction along $\B f$.

As an example consider the natural inclusion $\Lie\hookrightarrow\Ass$ and the induced map of the cobar-constructions $\B\Ass\rightarrow\B\Lie$. Recall that an $A_\infty$-algebra is an algebra over the operad $\B\Ass$ whereas a $C_\infty$-algebra is an algebra over $\B\Lie$. We deduce that a minimal model of a $C_\infty$-algebra may be calculated
using either sums over $\Lie$-decorated trees or sums over planar trees; this recovers a result of \cite{GeC}.
\end{rem}

\section{Homotopies for operad maps and equivalence of operadic algebras}\label{sechomotopy}
In this section we will set up the framework for studying the notion of a homotopy equivalence of operadic algebras. This is a natural generalization of the notion of a homotopy equivalence of $A_\infty, L_\infty$ and $C_\infty$-algebras. We will start with a general discussion of the Maurer-Cartan moduli space and the Sullivan homotopy; presumably this material is well-known to experts.

\subsection{Maurer-Cartan moduli space and the Sullivan homotopy}
Let $\g$ be a nilpotent dg Lie algebra or, more generally, a \emph{pro-nilpotent} Lie algebra, i.e. an inverse limit of nilpotent dg Lie algebras.
\begin{defi}
An element $x\in \g^1$ is called a \emph{Maurer-Cartan} element if it satisfies the following \emph{master equation}:
\[dx+\frac{1}{2}[x,x]=0.\] The set of Maurer-Cartan elements will be denoted by $\MC(\g)$.
\end{defi}
Let $G$ be the Lie group obtained by exponentiating the even part $\g^0$ of the Lie algebra $\g$.
Formally $G$ could be defined as the set of group-like elements in the universal enveloping algebra
$\widehat{U}(\g)$
of $\g$ completed with respect to its maximal ideal. Whenever we consider linear combinations of monomials in $\g$ or $G$ these will assumed to be taken in this completed universal enveloping algebra.

The group $G$ acts on $\MC(\g)$ by the formula
\[g(x)=gxg^{-1}-dg\cdot g^{-1}.\]
We say that two Maurer-Cartan elements are \emph{equivalent} if they lie in
the same $G$-orbit.

It is convenient to introduce the Lie algebra $\tilde{\g}$, the semi-direct product of $\g$ and a one-dimensional Lie algebra spanned by an odd symbol $d$. By definition for $a\in \g$ we have $[d,a]:=da$. For an element $x\in\g$ denote by $\tilde{x}$ the element $d+x\in\tilde{g}$. Then an odd element $x\in \g$ is Maurer-Cartan if and only if $[\tilde{x},\tilde{x}]=0$. Let $g\in G$ be viewed as an element in $\widehat{U}(\tilde{\g})$; since $d(g)=dg-gd$ we have $g dg^{-1}=d-dg\cdot g^{-1}$. It follows that the corresponding action of $G$ translates into the formula $g(\tilde{x})=g\tilde{x}g^{-1}$; in particular it is now obvious that this action is indeed well-defined.

Let $D:=\ground[z,dz]$, the differential graded commutative algebra generated by an even symbol $z$ and an odd symbol $dz$ with differential $d(z)=dz$ and $d(dz)=0$; this is just the polynomial de Rham algebra on the unit interval. Note that the specializations $z=0$ and $z=1$ determine two algebra maps $\ev_{0,1}:D\ra\ground$.

We say that two Maurer-Cartan elements $x_0$ and $x_1$ are  \emph{(Sullivan) homotopic} if there exists a Maurer-Cartan element $X\in\g\otimes D$ such that $(id\otimes \ev_0)(X)=x_0$ and $(id\otimes \ev_1)(X)=x_1$. The element $X$ is called a \emph{homotopy} between $x_0$ and $x_1$.

Let $X=x(z)+y(z)dz$ be a homotopy as defined above.
Here $x(z)\in\g^1[z]$ and $y(z)\in\g^0[z]$.
The following result is immediate from the definition.
\begin{lem} A homotopy between $x_0$ and $x_1$ is equivalent to the following system of equations.
\begin{eqnarray}\label{homotop1}[\tilde{x}(z),\tilde{x}(z)]=0.\\
\label{homotop2}\partial_z\tilde{x}(z)=[y(z),\tilde{x}(z)].\end{eqnarray}
(together with the boundary conditions $\tilde{x}(0)=\tilde{x}_0$ and
$\tilde{x}(1)=\tilde{x}_1$.)
\end{lem}
\begin{rem} The above identities could, of course, be rewritten as
\[dx(z)+\frac{1}{2}[x(z),x(z)]=0.\]
\[\partial_zx(z)=-dy(z)+[y(z),x(z)];\] however we find it more convenient to work with $\tilde{x}$ than with $x$.
\end{rem}
We see, therefore, that a homotopy is a one-parameter deformation $x(z)$ of a Maurer-Cartan element such that the differential equation $\partial_zx(z)=-dy+[y,x]$ holds. The following important result seems to be well-known although we are not aware of any published proof. It appears in an unpublished manuscript of Schlessinger and Stasheff \cite{SS}.
\begin{theorem}\label{MC}
Two Maurer-Cartan elements $x_0$ and $x_1$ are equivalent if and only if they are homotopic.
\end{theorem}
\begin{proof}
Let $x_0$ and $x_1$ be equivalent Maurer-Cartan elements. Then there exists an element $\xi\in\g$ such that $e^\xi \tilde{x}_0e^{-\xi}=\tilde{x}_1$. Set $\tilde{x}=e^{z\xi} \tilde{x}_0e^{-z\xi}$ and $y=\xi$; this clearly satisfies (\ref{homotop1}) and (\ref{homotop2}) and thus establishes a homotopy between $x_0$ and $x_1$.

Now suppose that $x_0$ and $x_1$ are homotopic, so that there exist
$\tilde{x}(z)\in\tilde{\g}[z]$ and $y(z)\in\g[z]$ for which
(\ref{homotop1}) and (\ref{homotop2}) hold. It is easy to check
that $\tilde{x}(z)$ is determined by $y(z)$ together with the
boundary condition $\tilde{x}(0)=\tilde{x}_0$.

If we could find an element $g(z)\in G[z]$ such that
$$\partial_zg(z)=y(z)g(z) \quad \text{and} \quad g(0)=1$$
then (\ref{homotop2}) would be satisfied with
$g(z)\tilde{x}_0 g(z)^{-1}$ in place of $\tilde{x}(z)$. By uniqueness
we could deduce
 $\tilde{x}=g(z)\tilde{x}_0 g(z)^{-1}$, and
$\tilde{x}_0$ and $\tilde{x}_1$ would be equivalent as desired.

We define $g(z)$ as a path ordered exponential:
\[g(z)\,=\,\operatorname{P}\exp\int_0^zy(t)dt\,:=\, 1+\sum_{n=1}^{\infty}\int_{0\leq t_1\leq \ldots \leq t_{n} \leq z}
y(t_n)\ldots y(t_1) \;\Pi dt_i.\]
Since $\g$ is pronilpotent we have $g(z)\in\widehat{U}(\g)[z]$, and the
differential equation $\partial_zg(z)=y(z)g(z)$ is clearly satisfied.
To complete the proof it remains to show that $g$ is group-like, i.e., that
$\Delta g=g\otimes g$; we follow the argument of Connes and Marcolli
\cite[Proposition 2.9]{CM}.
It suffices to prove
that the coefficient of $z^n$ in $\Delta g-g\otimes g$ is zero for all $n\geq 0$. We have
 \begin{eqnarray*}
\partial_z(\Delta g-g\otimes g) &=&
\Delta(\partial_z g) - \partial_z g\otimes g - g\otimes\partial_z g \\
&=& \Delta(y g) - yg\otimes g - g\otimes\ y g \\
&=& (y\otimes 1 + 1\otimes y) (\Delta g-g\otimes g).
\end{eqnarray*}
It follows that  $\partial_z^n(\Delta g-g\otimes g)$ is divisible by $\Delta g-g\otimes g$ for $n\geq 0$. It remains to observe that the constant term of $\Delta g-g\otimes g$ is zero because
$g(0)=1$.
\end{proof}
\begin{rem}
Our proof suggests that Theorem \ref{MC} could extend to  a not necessarily nilpotent dg Lie algebra supplied with a Banach norm. We plan to return to this issue in a future work.
\end{rem}

\subsection{Weak equivalence of operadic algebras}\label{equiva}
There is a natural notion of \emph{homotopy} between maps of operads; let $\O$
and $\O'$ be dg operads  and $f_{0},f_{1}:\O\ra\O'$ be two maps between them.
\begin{defi}
We say that $f_0$ and $f_1$ are \emph{homotopic} if
if there exists a map
$f:\O\ra \O'\otimes D$ of operads such that $f_0=(id\otimes \ev_{0})\circ f$ and $f_1=(id\otimes \ev_{1})\circ f$.
\end{defi}

Markl, Shnider and Stasheff call this relation `elementary homotopy', reserving the term `homotopy' for its transitive closure
\cite[Definition 3.121]{MSS}. They prove the following:

\begin{prop} \label{integrate1}
 Two homotopic maps $f_{0}, f_{1}:\O\ra\O'$ of operads induce the same maps
between the homology operads $H(\O)$ and $H(\O')$.
\end{prop}

\begin{proof}
Let $f=a(z)+b(z) dz$ be a homotopy
between $f_0$ and $f_1$. The
compatibility of $f$ with the differential implies
$(d\circ b)(z)+\partial_z a(z) = -(b\circ d)(z)$. Integrating and noting that $-\int_0^1 \partial_z a(z) = a(0)-a(1)=f_0-f_1$,
we find that $s=\int_0^1 b(z) dz$ is a chain homotopy
between $f_0$ and $f_1$. It follows that $f_0$ and $f_1$ induce the same homomorphism $H(\O_1)\ra H(\O_2)$ as claimed.
\end{proof}

\begin{defi}\label{homotopy}
Let $V$ be a dg vector space.
Two $\O$-algebras determined by operad maps $f_0:\O\ra\E(V)$ and $f_1:\O\ra\E(V)$ are called \emph{homotopy equivalent} if $f_0$ and $f_1$ are homotopic.
\end{defi}
The following is an immediate consequence of Proposition \ref{integrate1}.
\begin{cor}\label{integrate}
Two homotopy equivalent $\O$-algebra structures on a dg vector space $V$ give rise to the same $H(\O)$-algebra structure on $H(V)$.
\end{cor}

One could expect the above definition to give the correct notion of equivalence only for cofibrant admissible operads $\O$, i.e. those whose underlying operads of graded vector spaces are free. We will be interested in this notion in the special case when $\O=\B\P$ is the cobar-construction of an admissible operad $\P$. The notions of a $\B\P$-algebra and of an equivalence between two $\B\P$-algebras admits a reformulation in terms of a certain Maurer-Cartan moduli space as follows.

First of all,
 a (graded or super-) \emph{derivation} of a $\P$-algebra $A$ is an even map $f:A\rightarrow A$ such that for $p\in\P(n)$, \[f(p(a_1,\ldots,a_n))=\sum(-1)^{|f|(|p|+|a_1|+\ldots+|a_{i-1}|)}p(a_1,\ldots, a_{i-1},f(a_{i}),a_{i+1},\ldots,a_n).\] It is clear that a derivation of $A$ could be viewed as an infinitesimal automorphism of $A$ in agreement with the familiar notion. It is further clear that the space spanned by all graded derivations forms a Lie superalgebra.

Now let $V$ be a pro-finite vector space, i.e. an inverse limit of finite-dimensional vector spaces. For example, the linear dual to any (not necessarily finite-dimensional) vector space is pro-finite. Consider the pro-free $\P$-algebra on a dg vector space $ V$:
\[\hat{T}_\P(V)=\prod_{i=1}^\infty \P(i)\otimes_{\ground[\S_i]}( V)^{\hat{\otimes} i}.\] It is clear that $\hat{T}_\P(V)$ has a linear topology such that the structure maps $\P(n)\otimes[\hat{T}_\P(V)]^{\otimes n}\rightarrow \hat{T}_\P(V)$ are continuous. It is likewise clear that all
continuous derivations of $\hat{T}_\P(V)$ are determined by their values on $ V$. In other words, the space of all such derivations is isomorphic to $\prod_{i=1}^\infty\Hom(V,\P(i)\otimes_{\ground[\S_i]}( V)^{\hat{\otimes} i})$. The elements of the component $\Hom(V,\P(i)\otimes_{\ground[\Sigma_i]}( V)^{\hat{\otimes} i})$ will be called derivations of \emph{order i} for obvious reasons. The commutator of two derivations of orders $i$ and $j$ has order $i+j-1$.

Denote by $L(\hat{T}_\P(V))$ the space spanned by derivations of order $\geq 2$; it is a pro-nilpotent Lie superalgebra.

Then we have the following result
\begin{prop}\label{corresp}
Let $\P$ be an admissible operad. Then there is a one-to-one correspondence between the set of $\B\P$-algebra structures on a dg vector space $V$ and the set of Maurer-Cartan elements in $L(\hat{T}_\P(V^*))$. Furthermore, two Maurer-Cartan elements are equivalent if and only if the corresponding $\B\P$-algebras are homotopy equivalent.
\end{prop}
\begin{proof}
Since the operad $\B\P$ is freely generated by the collection $\{\Pi\overline{\P}(n)^*\}$  we see that an operad map $\B\P\rightarrow \E(V)$ determining a $\B\P$-algebra structure on $V$ is specified by a collection of $\S_n$-equivariant maps $\Pi\P(n)^*\rightarrow \Hom(V^{\hat{\otimes} n}, V),\, n=2,3,\ldots$.
Using the canonical isomorphism between $\S_n$-invariants and $\S_n$-coinvariants we deduce that the set of all operad maps $\B\P\rightarrow \E(V)$ (not assuming compatibility with the differential in $\B\P$)  is in one-to-one correspondence with the set of (possibly infinite) collections of odd maps $V^*\rightarrow \P(n)\otimes_{\ground[\Sigma_n]} (V^{\hat{\otimes} n})^*$, i.e., with odd elements in $L(\hat{T}_\P(V^*))$. Finally, compatibility with the differential is equivalent to the Maurer-Cartan identity $[\xi,\xi]=0$.

It follows immediately from definitions that two $\B\P$-algebra structures are homotopy equivalent if and only if the corresponding Maurer-Cartan elements are homotopic and thus, by Theorem \ref{MC} if and only if they are equivalent.
\end{proof}
\begin{rem}
Proposition \ref{corresp} shows that  our notion of homotopy equivalence between two $L_\infty, C_\infty$ or $A_\infty$-algebras (supported on the same dg vector space $V$) coincides with the usual notion of infinity-isomorphism, cf. for example \cite{keller}, \cite{HL1, HL2}.
\end{rem}
\begin{rem}
A version of the above proposition (without the statement about homotopy equivalence) appeared in \cite{getjon}, Proposition 2.15. The main difference with our approach is that in the cited reference the result is formulated in terms of the coalgebra (continuously) dual to the topological algebra $\hat{T}_\P(V^*)$. The approach taken here has been used in \cite{HL1, HL2} in the special cases of $L_\infty, C_\infty$ and $A_\infty$-algebras.
\end{rem}
\begin{cor}\label{equiv}
The relation `homotopy equivalence' on the set of all $\B\P$-algebra structures on a given dg vector space $V$ is an equivalence relation.
\end{cor}

\subsection{Minimal models and homotopy equivalence}In this subsection we prove that the minimal model of a $\B\P$-algebra is in fact homotopy equivalent to the original $\B\P$-algebra structure. In fact, we establish a slightly more general result which says, informally, that any two choices of a Hodge decomposition on an $\O$-algebra $V$ are homotopic as maps out of $\bv\O$.
\begin{theorem}\label{weakly}
Let $\O$ be an admissible dg operad,  and $f:\O\rightarrow \E(V)$ be a map of dg operads determining an $\O$-algebra structure on a dg vector space  $V$. Let $h_1, h_2:\BV\ra\E(V)$ be the maps of dg operads determined by the choice of two Hodge decompositions on $V$. Then the $\bv\O$-structures on $V$ corresponding to the operad maps $h_1\circ i$ and $h_2\circ i$ are homotopy equivalent.
\end{theorem}
\begin{proof}
The idea is that, while the operad $\bv\O$ is too complicated to construct the required homotopy directly, one could attempt to embed $\bv\O$ into a simpler operad, for example $\BV\O$ and construct the required homotopy as a map from that bigger operad.

It turns out that $\BV\O$ cannot do the job required.   We now construct another operad $\bigBV\O$ such that $\bv\O\hookrightarrow\bigBV\O\rightarrow\BV\O$. Informally speaking, $\bigBV\O$ is `sufficiently cofibrant' so that the required homotopy exists while it is simple enough for the homotopy to be written down explicitly.

The new operad $\bigBV\O$ is generated by the operad $\O$ together with
free noncommuting generators $s$ and $t$ and with differential given by $d(s) = \mathbf{1}-t^2$ and $d(t)=0$.
There is map of operads $\bigBV\O\rightarrow\BV\O$ given by imposing the additional relations $s^2=0$, $t^2=t$ and $st=ts$.

Furthermore, the operad $\bv\O$ is a suboperad of $\bigBV\O$. To see that note that the latter consists of $\O$-decorated reduced trees, with internal edges labelled by nonempty words in $s$ and $t$, and extremities labelled by arbitrary words in $s$ and $t$. Then $\bv\O$ is identified with the suboperad of $\bigBV\O$ where internal edges are labelled by $s$ or $t^2$, and extremities are labelled by $t$.

Let $S$ and $T$ be the two operators on $V$ determined by a given Hodge decomposition. Denote by $\hat{f}:\bigBV\O\rightarrow\E(V)$ the operad map which equals $f$ when restricted to $\O$ and for which $\hat{f}(s)=S$, $\hat{f}(t)=T$. Similarly denote by $\hat{g}:\bigBV\O\rightarrow\E(V)$ the operad map which again equals $f$ when restricted to $\O$ and for which $\hat{g}(s)=0$, $\hat{g}=\text{id}$.

Then the following formulas determine a homotopy.
\[u:\bigBV\O\rightarrow\E(V)\otimes\ground[z,dz]:\]
\[u(s) = S(1-z^2);\]
\[u(t)= T +  (1-T) z - S dz.\]
Since $\bigBV\O$ is free over $\O$ one should only check that these formulas are compatible with differentials in $\bigBV\O$ and $\E(V)\otimes\ground[z,dz]$ which is straightforward. Moreover setting $z=0$ and $z=1$ we recover the maps $\hat{f}$ and $\hat{g}$ respectively. It follows that the maps $\hat{f}$ and $\hat{g}$ are homotopic as required. Thus their restrictions to $\bv\O$ are homotopic, and the desired
result is a consequence of Corollary~\ref{equiv}.
\end{proof}
\begin{rem}
Note that once we do not require $t$ to be an idempotent in $\bigBV\O$, we must use the equation $d(s) =\mathbf{1}-t^2$ instead of $d(s)=\mathbf{1}-t$. This is because when we glue together two trees at extremities both labelled by $t$, the new internal edge is now labelled by
$t^2$, not $t$.
\end{rem}
\begin{rem}
The proof of Theorem \ref{weakly} consisted in construcing a homotopy between two Hodge decompositions. The homotopy, however, passes through structures which are 
 not themselves Hodge decompositions, since no relations are imposed on the generators $s$ and $t$ of $\bigBV\O$. We do not know a good interpretation of such `generalized Hodge decompositions'.
\end{rem}

\begin{cor}
Let $\O$ be an admissible dg operad,  and $f:\O\rightarrow \E(V)$ be a map of dg operads determining an $\O$-algebra structure on a dg vector space  $V$. Let $\tilde{f}:\O\rightarrow \E(V)$ be the minimal $\O$-algebra structure on $V$ associated with a given splitting $k:\bv\ra\O$ of the canonical resolutions $\bv\O\ra\O$ and a canonical Hodge decomposition of $V$. Then the two $\O$-algebra structures on $V$ corresponding to $f$ and $\tilde{f}$ are homotopy equivalent. In particular, any two minimal models are likewise homotopy equivalent.
\end{cor}
\begin{proof}
Let $\epsilon_1:\bv\O\ra \E(V)$ be the operad map associated with the given canonical Hodge decomposition on $V$. Let $\epsilon_2:\bv\O\ra\E(V)$ be the map associated with the \emph{trivial} Hodge decomposition of $V$. That means that $\epsilon_2=\epsilon\circ i$ where $i:\bv\O\hookrightarrow\BV\O$ is the canonical embedding and $\epsilon:\BV\O\ra\E(V)$ is given by $\epsilon(s)=0$ and $\epsilon(t)=\operatorname{id}$. By Theorem \ref{weakly} the maps $\epsilon_1$ and $\epsilon_2$ are Sullivan homotopic so the desired statement follows.
\end{proof}
\begin{rem}
Our result is more general than that usually called the `minimal model theorem' in the literature; it specializes to the so-called `decomposition theorem' in the case of $A_\infty$-algebras, cf. \cite{Kajiura, Lef}. One can rephrase it by saying that for any operad $\O$ for which there is a splitting of the canonical map $\bv\O\rightarrow \O$ (e.g. $\O$ could be a cobar-construction of an admissible operad), any $\O$-algebra is infinity-isomorphic (=homotopy equivalent)
to an $\O$-algebra of a special type, namely a direct sum of an $\O$-algebra with vanishing differential and what was called in \cite{Kajiura} a \emph{linear contractible } $\O$-algebra.
 From this result it is not hard to deduce that different minimal models of an $\O$-algebra $V$ \emph{understood as $\O$-algebra structures on $H(V)$}
 are homotopy equivalent.
\end{rem}

\begin{rem} Let $V$ be an algebra over an admissible
operad $\O$. Then
the minimal model, viewed as an $\O$-algebra structure on
$H(V)$, is a lift of the $H(\O)$-algebra structure on $H(V)$ induced by the original $\O$-algebra $V$. Indeed, the $\O$-algebra structures on $V$ are homotopy equivalent by the preceding corollary and hence coincide in homology, by
Proposition~\ref{integrate}.
\end{rem}

\section{Minimal models of algebras over modular operads}\label{modular}
In this section we construct minimal models for algebras over modular operads and prove that these are unique up to a non-canonical isomorphism. We restrict ourselves to giving the relevant definitions and formulations; the proofs will be omitted as they are completely parallel to those in the non-modular case.

We refer the reader to \cite{GeK} for generalities on modular operads and \cite{CL} for the notion of the $\BV$-resolution of algebras over modular operads; we shall liberally use terminology and notation from these two sources. For a dg vector space $V$ with a symmetric inner product $\langle, \rangle$ of
even degree we will still denote by $\E(V)$ the modular operad of endomorphisms of $V$, with $\E(V)((g,n)):=
V^{\otimes n}\cong \Hom(V^{\otimes n-1}, V)$.   The self-glueing maps in $\E(V)$ are determined by the inner product in $V$.

For a modular operad $\O$ the structure of an algebra over $\O$ on $V$ is a map of  modular operads $\O\ra\E(V)$. A \emph{Hodge decomposition} of an algebra $V$ over $\O$ is a pair of operators $s$ and $t$ as in Proposition \ref{dualalg1} compatible with the inner product:
\[\langle s(a),b\rangle=(-1)^{|a|}\langle a,s(b)\rangle;\]
\[\langle t(a),b\rangle=\langle a,t(b)\rangle.\] As before, a Hodge decomposition will be called \emph{canonical}
if $dt=0$
and \emph{trivial}
if $t=\id_V$.
 It is easy to see that a canonical Hodge decomposition always exists; indeed, following the
procedure in Example~\ref{Hodge} we find that $U$ is a maximal
isotropic subspace of $W^\perp$ and just need to choose $U'$ to be an isotropic
complement. The diagram (\ref{operadmaps}) continues to hold in the modular context.

\begin{rem} In order to handle inner products of \emph{odd} degree (and to
discuss Feynman transforms) we are obliged to use the language of twisted modular operads. Let $V$ be a dg vector space with a
symmetric inner product of even or odd degree. For $d=0,1$, we define the twisted endomorphism operad $\E_{d}(V)$ with components $\E_{d}(V)((g,n))=\Pi^d(\Pi^d V)^{\otimes n}$. If the inner product on
$V$ has degree $d+d'$ then
$\E_{d}(V)((g,n))
\cong \Pi^{d'}\Hom((\Pi^{d'}V)^{\otimes {n-1}},\Pi^{d'}V)$, and
$\E_{d}(V)$ is a modular $\De\otimes{\mathfrak K}^{d'}$-operad, where
$\De$ is the determinant cocycle $\Det^d(G)=\Det^d(H_1(G))$ and $\mathfrak K$ is the dualizing cocycle  cf. \cite{GeK} and \cite{bar}.

So given an arbitrary modular $\De\otimes{\mathfrak K}^{d'}$-operad $\O$,
we define an $\O$-algebra to be a dg space $V$ equipped with a symmetric inner product of degree $d+d'$ together with a map $\O\rightarrow\E_d(V)$.

 To alleviate the notation we will suppress any explicit mentioning of twisting whenever possible if the context  allows one to reconstruct it; for example  for a cocycle $\mathfrak D$ and a modular $\mathfrak D$-operad $\O$ we will write $\F\O$ instead of $\F_{\mathfrak D}\O$ to denote the Feynman transform
 of $\O$.
\end{rem}

If $V$ has vanishing differential then the corresponding $\O$-algebra structure is called \emph{minimal}. If the canonical operad map $\bv\O\ra\O$ admits a splitting then we can construct a minimal model of $V$ as in
Section~\ref{main}.

So let $V$ be an algebra over a modular operad $\O$ with a fixed
canonical Hodge decomposition of $V$. Suppose that $\O$ is the Feynman transform $\F\P$ of a modular $\Det^d$-operad $\P$. Recall
that as a graded modular operad, i.e. forgetting the differential, $\F\P$ is
free over the stable $\S$-module $\P^*$;
 therefore the operad map $\F\P\rightarrow \E_{d}(H(V))$ providing a minimal model of $V$ is determined by a collection of maps \[\Pi\P((g,n))^*\rightarrow \Pi\E_{d}(H(V))((g,n))=\Hom(H(\Pi V)^{\otimes {n-1}},H(\Pi V)).\]  These maps will be called the structure maps of the corresponding minimal model.

\begin{example}\
\begin{enumerate}
\item Let $\P=\underline{\Comm}$, the trivial modular extension of the operad $\Comm$, so that
\[
\P((g,n))=\P((g,n))^*=\begin{cases}{\ground \text{~if~} g=0}\\{0 \text{~if~} g\neq 0}\end{cases}.
\]
Then we  obtain a collection of $\S_n$-equivariant maps $m_n:\Hom(H(\Pi V)^{\otimes n},H(\Pi V))$ which determine the structure of a minimal symplectic (or cyclic) $L_\infty$-algebra on $H(V)$.
\item Let
$\P=\overline{\Comm}$, modular closure of the operad $\Comm$, so that $\P((g,n))=\P((g,n))^*=\ground$. There results a collection of $\S_n$-equivariant maps $m_{g,n}:\Hom(H(\Pi V)^{\otimes n},H(\Pi V))$ which determine the structure of a minimal quantum $L_\infty$-algebra on $H(V)$, also known as a loop homotopy algebra, cf. \cite{Markl}. This structure first arose in closed string field theory \cite{Zwi}.
\item
Let $\P=\underline{\Ass}$, the trivial modular extension of the operad $\Ass$, so that
\[
\P(g,n)=\begin{cases}{\ground[\S_n/\mathbb{Z}_n]\text{~if~} g=0}\\{0 \text{~if~} g\neq 0}\end{cases}.
\]
We  obtain a collection  of maps $m_n:\Hom(H(\Pi V)^{\otimes n},H(\Pi V))$ which determine the structure of a minimal symplectic (or cyclic) $A_\infty$-algebra on $H(V)$; this structure was introduced by Kontsevich in \cite{K1,K2}.
\item
Let $\P=\overline{\Ass}$, the  modular closure of the operad $\Ass$, so that $\P((g,n))=\ground[\S_n/\mathbb{Z}_n]$.  There results a collection of maps $m_n:\Hom(H(\Pi V)^{\otimes n},H(\Pi V))$ which determine a structure  on $H(V)$ which is natural to call a \emph{quantum $A_\infty$-algebra}. More information on other modular extensions of $\Ass$ and their connections with various compactifications of moduli of Riemann surfaces could be found in \cite{CL}.
\end{enumerate}
\end{example}
Let $G$ be a stable graph (see \cite{GK}); for a vertex $v\in\Ve(G)$ the set of half-edges around $v$ is denoted by $\Flag(v)$.
For a stable $\S$-module $\V$, we denote by $\V((\G))$ the space of $\V$-decorations on $G$; in other words, $G((\V))=\otimes_{v\in\Ve{G}}\V((\Flag(v)))$.

Let $\P$ be a modular $\De$-operad.
For any stable graph $G$ with $n$ legs, $\P$ determines a homomorphism $\mu_G:{\De}((G))\otimes\P((G))\ra\P((n))$ which corresponds to taking operadic compositions in $\P((G))$ along the internal edges of $G$; the dual map will be denoted by $\mu_G^*$.

Let $V$ be a $\F\P$-algebra. Note that $V$ has an inner product of
degree $1-d$. Given any twisted $\P^*$-decoration
$x\in \De((G))\otimes\P^*((G))$ on a stable graph $G$ with $n$ legs, we define
the BV-amplitude $\zeta_V(x)\in\E_d(V)((n))$ as follows. We replace the
$\P^*$-decorations at the vertices by
$\E_d(V)$-decorations via the action of
$\P^*\subset\F\P$ on $V$. So each vertex of $G$ of valence $m$
 has attached to it
a tensor in $\Pi^{d(1+m)}V^{\otimes m}$. Then we assign
`propagators' $\langle ? ,s(?)\rangle \in \Pi^{d}(V^*\otimes V^*)$ to each internal edge and $t\in V\otimes V^*$ to each leg, and contract along all edges (including legs) to obtain a tensor $\zeta_V(x)\in\Pi^{d(1+n)}V^{\otimes n}$.


We have the following result which is a modular analogue of Theorem \ref{main1}.
Its proof is also an exact analogue, making use of the notions of
$BV$-resolutions of modular operads and $BV$-graphs developed in
\cite{CL}.

\begin{theorem}\label{main2}Let $\P$ be a modular $\Det^d$-operad
Then any choice of a canonical Hodge decomposition on an $\F\P$-algebra $V$  gives rise to a $\F\P$-algebra structure on $H(V)$ that is a minimal model of $V$.
The structure maps of this minimal model are
$\sum_G \zeta_V\circ \mu_G^*:
\P((g,n))^*\rightarrow \E_d(V)((g,n))$, where the sum is extended over
all stable $n$-graphs $G$ of genus $g$.
\end{theorem}

\begin{rem}
Let $\C$ be a cyclic operad and $\P=\underline{\C}$ its trivial
modular extension, so that all `self-composition' maps
in $\P$ are zero. Then the sum in Theorem~\ref{main1} may
be restricted to trees, so we recover the formula of
Theorem~\ref{main}. This relationship between the minimal model construction in the operad and modular operad cases stems from the fact that the Feynman transform $\F\P$ is (up to a twist) the modular closure of the cyclic operad $\B\C$.

Taking for example $\C=\Ass$, we see that the explicit minimal model for an symplectic $A_\infty$ algebra is calculated in precisely the same way as for an $A_\infty$-algebra, once
a Hodge decomposition of $V$ compatible with inner product is chosen, c.f. \cite{Kajiura}.
\end{rem}

\subsection{Weak equivalence of algebras over modular operads}
Definition \ref{homotopy} of the homotopy equivalence of $\O$-algebra structures on $V$ carries over to the modular context verbatim.
One would, of course, expect it to give the correct notion of equivalence only for `cofibrant' modular operads $\O$, i.e. those whose underlying modular operads of graded vector spaces are free. We will be interested in this notion in the special case when $\O=F\P$ is the Feynman transform of a  modular operad $\P$. The notions of an $\F\P$-algebra and of an equivalence between two $\F\P$-algebras admits a reformulation in terms of a certain Maurer-Cartan moduli space as follows.

 Consider the following pro-finite vector space:
\[L_\P(V):=\prod_{g,n=0}^\infty [\P((g,n))\otimes V^{\hat\otimes n}]^{\S_n}.\]
This space is a modular analogue of $L(\hat{T}_\P(V))$  considered in Section \ref{equiva}. It has the structure of a pronilpotent dg Lie algebra, cf. \cite{bar}. We have the following modular analogue of Proposition \ref{corresp}.
\begin{prop}\label{corresp1}
Let $\P$ be a modular operad. Then there is a one-to-one correspondence between the set of $\F\P$-algebra structures on a dg vector space $V$ and the set of Maurer-Cartan elements in $L_\P(V)$. Furthermore, two Maurer-Cartan elements are equivalent if and only if the corresponding $\F\P$-algebras are homotopy equivalent.
\end{prop}
\begin{proof}
The first statement of the proposition is essentially Theorem 1 of \cite{bar} whereas the second one is follows immediately from definitions.
\end{proof}
\begin{rem}As in the non-modular case Proposition \ref{corresp1} implies that the homotopy equivalence between algebras over modular operads of the form $\F\O$ is an equivalence relation.
\end{rem}

\begin{theorem}\label{weaklymod}
Let $\O$ be an modular operad,  and $f:\O\rightarrow \E(V)$ be a map of modular operads determining a structure of an $\O$-algebra on a dg vector space  $V$ with an inner product. Let $h_1, h_2:\BV\O\ra\E(V)$ be the maps of modular operads determined by the choice of two different Hodge decompositions on $V$. Then the $\F\F\O$-structures on $V$ corresponding to the operad maps $h_1\circ i$ and $h_2\circ i$ are homotopy equivalent.
\end{theorem}

\begin{cor}
Let $\O$ be an modular operad,  and $f:\O\rightarrow \E(V)$ be a map of modular operads determining an $\O$-algebra structure on a vector space  $V$ with an inner product. Let $\tilde{f}:\O\rightarrow \E(V)$ be the minimal $\O$-algebra structure on $V$ associated with a given splitting $k:\O\ra\F\F\O$ of the canonical resolutions $\F\F\O\ra\O$ and a canonical Hodge decomposition of $V$. Then the two $\O$-algebra structures on $V$ corresponding to $f$ and $\tilde{f}$ are homotopy equivalent. In particular, any two minimal models are likewise homotopy equivalent.
\end{cor}

\section{Gauge independence of Kontsevich's dual construction}
Recall the operadic formulation of Kontsevich's dual construction \cite{CL}. Let $\O$ be a modular operad, $V$ be a \emph{contractible} dg vector space with an inner product and a structure of an $\O$-algebra. The choice of a canonical Hodge decomposition (which amounts to the choice of a contracting homotopy $s$ such that $s^2=0$ in this case) determines the structure of a $\F^\vee\O$-algebra on $V$. Here $\F^\vee\O$ is the dual Feynman transform of $\O$; to form it we first add a unit $\bf 1$ to $\O$ to obtain an extended operad $\O_+$ and then freely adjoin an odd element $s\in \F^\vee\O((0,2))$ subject to the relation $s^2=0$. We can thus write $\F^\vee\O=\O[s]_+/s^2$; the differential in $\F^\vee\O$ extends the differential on $\O$ and $d(s)=\bf 1$.

The main property of $\F^\vee\O$ is that there is a linear duality between $\F^\vee\O((g,0))$ and $\F\O((g,0))$. Taking the Feynman amplitude $G\mapsto Z_{V}^{\F^\vee\O}(G)$
determines a cocycle on $\F^\vee\O((0))$ (or equivalently, a cycle on $\F\O((0))$). We will denote the corresponding (co)homology class by $[V]$.

\begin{prop}\label{gauge}
The class $[V]$ corresponding to a contractible $\O$-algebra $V$ does not depend on the choice of the contracting homotopy $s$.
\end{prop}
\begin{lem}\label{sur}
Let $\tilde{\F}^\vee\O$ be the the extended modular operad $\tilde{F}^\vee\O=\O[s]_+$; the differential is specified by the same formula as in $\F^\vee\O$. Then the canonical map $\pi:\tilde{\F}^\vee\O\ra \F^\vee\O$ induces a surjective homomorphism
$H(\tilde{\F}^\vee\O((0)))\ra H(\F^\vee\O((0)))$.
\end{lem}
\begin{proof}
We have an isomorphism $\tilde{\F}^\vee\O((0))=\oplus_G \O\{G\}$, where the sum is
taken over all extended stable graphs with no legs and
$\O\{G\}=(\mathfrak{K}(G)\otimes\O((G)))_{\Aut(G)}$ is the
space of $\Aut(G)$-coinvariants of the (twisted)
space of $\O$-decorations on $G$; the differential is given
by contracting edges with the operad composition.
The dual Feynman transform $\F^\vee\O$ has an analogous decomposition where the sum is over stable graphs. The
restriction of the canonical map $\pi$ to $\O\{G\}$ is the identity
for stable $G$ and $0$ otherwise. While the obvious splitting
map $\F^\vee\O((0))\ra \tilde{\F}^\vee\O((0))$ is not a map
of operads, it at least commutes with differentials.
\end{proof}

\begin{rem}
In fact
$\pi$ induces an isomorphism
$H(\tilde{\F}^\vee\O((g,n)))\cong H(\F^\vee\O((g,n)))$ for
$(g,n)\neq(1,0)$. We do not require this stronger statement
for the proof below.
\end{rem}

\begin{proof}[Proof of Proposition \ref{gauge}]
Let $S$ and $S^\prime$ be two choices of a contracting homotopy in $V$. These determine two maps of modular operads $\tilde{\F}^\vee\O\ra \E(V)$ where $f(s)=S$ and $g(s)=S^\prime$. If we could prove that $f$ and $g$ are Sullivan homotopic then the desired statement would follow. In fact $f$ and $g$ are not homotopic due to the fact that $\F^\vee\O$ is not `cofibrant'; however the maps $f\circ\pi$ and $g\circ\pi$ are homotopic through the map $h:\tilde{\F}^\vee\O\ra \E(V)\otimes D$ specified by the formula
$h(s)=S+(S'-S)(z-Sdz)$. It follows by Proposition~\ref{integrate1} that $f\circ \pi$ and $g\circ \pi$ induce the same cohomology class on $\tilde{\F}^\vee\O((0))$, and an application of Lemma \ref{sur} finishes the proof.
\end{proof}

Having confirmed that the dual construction is gauge invariant, we can proceed
to show that it is also a Sullivan homotopy invariant of contractible algebras.

\begin{prop}
Let $\O$ be a modular operad, and let $V$ be a contractible dg vector space. Then Sullivan homotopic $\O$-algebra structures on $V$ give rise via the dual construction to the same class $[V]$ in $\F\O((0))$.
\end{prop}

\begin{proof}
Let $f:\O\to \E(V)\otimes D$ be a homotopy between the given $\O$-algebra
structures on $V$. Then a choice of contracting homotopy $s$ for $V$ allows us to extend $f$ to a homotopy $\tilde{f}:\F^\vee\O\to \E(V)\otimes D$ by taking
$\tilde{f}(s)=s\otimes 1$. The desired result is now a consequence of
Proposition~\ref{integrate1}
\end{proof}


\begin{thebibliography}{24}
%
\bibitem{AF} P. Aspinwall, L. M. Fidkowski. \emph{Superpotentials for quiver gauge theories.}  J. High Energy Phys.  2006,  no. 10, 047, 
\texttt{arXiv:hep-th/0506041}.

%
\bibitem{AK} P. Aspinwall, S. Katz. \emph{Computation of Superpotentials for D-Branes.}  Communications in Mathematical Physics 264, (2006) 227-253,
\texttt{arXiv:hep-th/0412209}.

%
\bibitem{bar} S. Barannikov. \emph{Modular operads and Batalin-Vilkovisky geometry.}
International Mathematics Research Notices, Vol. 2007, Article ID rnm075, 31 pages.
%
\bibitem{BV} J. M. Boardman, R.M. Vogt. \emph{Homotopy invariant algebraic structures on topological. spaces}. Springer Lecture Notes in Math. 347, Springer-Verlag, 1973.

%
\bibitem{CL} J. Chuang, A. Lazarev. \emph{Dual Feynman transform for modular operads},
Communications in Number Theory and Physics, to appear, \texttt{arXiv:0704.2561}.
%
\bibitem{GeC} X. C. Cheng, E. Getzler \emph{Transferring homotopy commutative algebraic structures}. \texttt{arXiv:math.AT/0610912}.
%
\bibitem{CM} A. Connes, M. Marcolli. \emph{Renormalization, the Riemann-Hilbert correspondence, and motivic
   Galois theory}, Frontiers in number theory, physics, and geometry. II.,
   Springer (2007), 617--713.

\bibitem{GeK} E.~Getzler, M.~M. Kapranov. \emph{Modular
    operads}. Compositio Math. 110 (1998), no.\ 1, 65--126.
%
\bibitem{GK} V. Ginzburg, M.~Kapranov. \emph{Koszul duality for
    operads}. Duke Math. J. 76 (1994), no.\ 1, 203-272.
%
\bibitem{getjon}
E. Getzler, J.D.S. Jones, \emph{Operads, Homotopy Algebra, and Iterated Integrals for double Loop Spaces.} \texttt{arXiv:hep-th/9403055}.
%
\bibitem{HL1}
A. Hamilton, A. Lazarev, \emph{Cohomology theories for homotopy algebras and noncommutative geometry}. \texttt{arXiv:0707.3937}
%
\bibitem{HL2} A. Hamilton, A. Lazarev, \emph{Symplectic $C_\infty$-algebras}. Moscow Mathematical Journal, to appear. \texttt{arXiv:0707.3951}
%
\bibitem{Kad}
T. V. Kadeishvili, \emph{The algebraic structure in the homology of an $A(\infty)$-algebra.} (Russian)  Soobshch. Akad. Nauk Gruzin. SSR 108 (1982), no. 2, 249-252, 1983.
%
\bibitem{Kajiura} H. Kajiura, \emph{Noncommutative homotopy algebras associated with open strings.}  Rev.Math.Phys. 19, 2007, 1-99.
\texttt{arXiv:math.QA/0306332}.
%
\bibitem{keller}
B. Keller, \emph{Introduction to $A$-infinity Algebras and Modules.} Homology, Homotopy and Applications, Vol. 3, 2001, No. 1, pp. 1-35.
%
\bibitem{K1}
M. Kontsevich. \emph{Feynman diagrams and low-dimensional topology.} First European Congress of Mathematics, Vol 2 (Paris, 1992), 97-121, Progr. Math., 120, Birkh\"auser, Basel, 1994.

\bibitem{K2}
M. Kontsevich. \emph{Formal noncommutative symplectic geometry.} The Gelfand Mathematical Seminars, 1990-1992, pp. 173--187, Birkh\"auser Boston, Boston, MA, 1993.
%
\bibitem{KS}M. Kontsevich, Y. Soibelman. \emph{ Deformations of algebras
over operads and the Deligne conjecture.} Conf\'erence Mosh\'e
Flato 1999, Vol. I (Dijon), 255-307, Math. Phys. Stud., 21, Kluwer
Acad. Publ., Dordrecht, 2000.
%
\bibitem{La} D. W. Barnes, L. A. Lambe.
\emph{A fixed point approach to homological perturbation theory.}
Proc. Amer. Math. Soc., 112(3), 881-892, 1991.
%
\bibitem{Laz} C. Lazaroiu. \emph{String field theory and brane superpotentials.}
J. High Energy Phys.
0110 (2001) 018,
 \texttt{arXiv:hep-th/0107162}.
 %
\bibitem{Lef}
K. Lef\`evre-Hasegawa, \emph{Sur les A-infini cat\'egories}, \texttt{arXiv:math.CT/0310337}.
%
\bibitem{Markl} M. Markl. \emph{Loop homotopy algebras in closed field theory}, Comm. Math. Phys., 221, 367-384, 2001.
%
\bibitem{markl} M. Markl. \emph{ Transferring $A_\infty$ (strongly homotopy associative) structures}, 
 Proceedings of the 25th Winter School Geometry and Physics. Palermo: Universita degli Studi di Palermo, 2006. Supplemento ai Rendiconti del Circolo Matematico,
 2006, 139-151.

%
\bibitem{MSS} M. Markl,  S. Shnider,  J. Stasheff. \emph{Operads in algebra, topology and physics},
Mathematical Surveys and Monographs, 96. American Mathematical Society, Providence, RI,  2002.
%
\bibitem{Mer} S. Merkulov. {\em Strongly homotopy algebras of a K\"ahler manifold}. Internat. Math. Res. Notices 1999, no. 3,
153-164.
%
\bibitem{SS} M. Schlessiger, J. Stasheff. \emph{Deformation theory and rational homotopy type}, unpublished manuscript.
%
\bibitem{Su} D. Sullivan. \emph{Infinitesimal computations in topology.}
Inst. Hautes \'Etudes Sci. Publ. Math. No. 47 (1977), 269-331 (1978).
%
\bibitem{To} A. Tomasiello. \emph{A-infinity structure and superpotentials.} J. High Energy Phys.
0109 (2001) 030, \texttt{arXiv:hepth/0107195}.
%
\bibitem{Zwi} B. Zwiebach. \emph{Closed string field theory: quantum action and the Batalin-Vilkovisky master equation.}
Nuclear Phys. B 390 (1993), no. 1, 33-152.
\end{thebibliography}
\end{document}